\documentclass[11pt]{nyjm}
\usepackage{amsfonts,amssymb,amsmath,amsthm,graphicx,amscd,float}
\usepackage{epic,eepic,epsfig,amsbsy,amsmath,color}
\usepackage[all]{xy}
\usepackage{texdraw}
\usepackage{url}

\usepackage{mathrsfs}
\usepackage{standalone}

\usepackage{tikz}
\usetikzlibrary{cd,positioning,knots,patterns,hobby}
\usepackage[normalem]{ulem}
\usepackage{comment}

    \newtheorem{theorem}{Theorem}[section]
    \newtheorem{prop}[theorem]{Proposition}
    \newtheorem{lemma}[theorem]{Lemma}
    
    \newtheorem{rem}[theorem]{Remark}
    \newtheorem{definition}[theorem]{Definition}
    \newtheorem{cor}[theorem]{Corollary}
    
    \newtheorem{fact}[theorem]{Fact}
    \newtheorem{scholium}{Scholium}

\newcommand{\ov}{\overline}
\newcommand{\cD}{\mathcal{D}}
\newcommand{\cS}{\mathcal{S}}

\newcommand{\bC}{\mathbb{C}}
\newcommand{\bZ}{\mathbb{Z}}
\newcommand{\BZ}{\mathbb{Z}}
\newcommand{\mZ}{\mathbb{Z}}

\newcommand{\cA}{\mathcal{A}}
\newcommand{\bi}{{\bf i }}
\newcommand{\be}{\begin{equation}}
\newcommand{\ee}{\end{equation}}
\newcommand{\inc}{\operatorname{inc}}

\newcommand{\im}{\operatorname{im}}

\newcommand{\Tr}{\operatorname{Tr}}

\newcommand{\cR}{\mathcal{R}}

\newcommand{\G}{SL_2(\mathbb{C})}
\newcommand{\PG}{PSL_2(\mathbb{C})}

\newcommand{\cL}{\mathcal{L}}

\newcommand\lk{\mathrm{lk}}

  \newcommand{\lcr}{\raisebox{-5pt}{\mbox{}\hspace{1pt}
                 \includegraphics{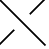}\hspace{1pt}\mbox{}}}
\newcommand{\ift}{\raisebox{-5pt}{\mbox{}\hspace{1pt}
                 \includegraphics{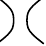}\hspace{1pt}\mbox{}}}
\newcommand{\zer}{\raisebox{-5pt}{\mbox{}\hspace{1pt}
                 \includegraphics{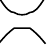}\hspace{1pt}\mbox{}}}

\def\al{\alpha}

\title[Skein Algebras at $4$th Roots of Unity]{Skein Algebras of Three-Manifolds at $4$th Roots of Unity}
\author[Frohman]{Charles Frohman}

\address{Department of Mathematics, The University of Iowa, Iowa City, IA
52242, USA}

\email{\tt charles-frohman@uiowa.edu}

\author[Kania-Bartoszynska]{Joanna Kania-Bartoszynska}

\address{National Science Foundation, Arlington, VA, 22230, USA} 

\email{\tt jkaniaba@nsf.gov}

\author[l\^{e}]{ Thang L\^{e}}

\address{Georgia Tech}

\email{\tt letu@math.gatech.edu}

\thanks{This material is based upon work supported by and while serving at the National Science Foundation. Any opinion, findings, and conclusions or recommendations expressed in this material are those of the authors and do not necessarily reflect the views of the National Science Foundation. }

\begin{document}
\maketitle
\begin{abstract}
This paper introduces an algebra structure on the part of the skein module of an arbitrary $3$-manifold $M$ spanned by links that represent $0$ in $H_1(M;\bZ_2)$ when the value of the parameter used in the Kauffman bracket skein relation is equal to $\pm {\bf i}$. It is proved that if $M$ has no $2$-torsion in $H_1(M;\bZ)$  then those algebras, $K_{\pm {\bf i}}^0(M)$, are naturally isomorphic to the corresponding algebras when the value of the parameter is $\pm 1$. This implies that  the algebra $K_{\pm{\bf i}}^0(M)$  is the unreduced coordinate ring of the variety of $P\G$-characters of $\pi_1(M)$ that lift to $\G$-representations. 
\end{abstract}

{\bf MSC2020 Classfication}: 57K31
\section{Introduction}
The skein module of a three-manifold is defined to be the quotient of the complex vector space with basis the  isotopy classes of framed links in the manifold by the sub-vector space spanned by  skein relations. The skein module admits an algebra structure for manifolds that are a product of a surface with an interval and for certain values of parameters used in the skein relation,  possibly with admissibility conditions on the links used.
Skein modules and skein algebras are central to quantum topology as a consequence of their relationship to the character variety. Namely, the Kauffman bracket skein algebra of a compact oriented three-manifold  $K_{-1}(M)$, when the variable in the Kauffman skein relation is set to $-1$, is the unreduced coordinate ring of the $\G$-character variety of the fundamental group of $M$ \cite{Bu,PS}. 

The Kauffman bracket skein module of any three-manifold is graded by $H_1(M;\bZ_2)$, 
\be K_\zeta(M)=\bigoplus_{x\in H_1(M;\bZ_2)}K_\zeta^x(M).\ee 
This paper introduces an algebra structure  on $K_{\pm {\bf i}}^0(M)$ for an arbitrary oriented three-manifold $M$ in {\bf Proposition} \ref{algstru}. The product is given in terms of ${\mathbb {Z}}_2$-linking numbers, by the distributive extension of the formula
 \be(\alpha,\beta) \mapsto (-1)^{\lk_2(L,L')} [L \sqcup L'].  \ee
 Here $\alpha$ and $\beta$ are skeins represented by disjoint links $L$ and $L'$ respectively.

A finite type surface is defined to be the complement of finitely many points in a closed, oriented surface. In a 2011 paper \cite{Ma}, March\'{e} proved that when the variable is set to $-{\bf i}$, the Kauffman bracket skein algebra of a finite type surface, $K_{-{\bf i}}(F)$,  is isomorphic to a twisted version of $K_{-1}(F)$. Specifically he constructs a noncommutative algebra $\cA$ isomorphic to  ${\mathbb{C}}H_1(F;\mZ_2)$ as a vector space. The algebra $\cA$ is graded by $H_1(F;\bZ_2)$.  He  constructs an isomorphism 
\be\label{march} \phi: K_{-{\bf i}}(F)\rightarrow (K_{-1}(F)\otimes \cA)_0 \ee where the subscripted $0$ denotes the $0$-graded part.

{\bf Theorem} \ref{bige} extends March\'{e}'s formula  to an isomorphism
\be \phi:K_{{\bf i}\zeta}(F) \rightarrow \left(K_\zeta(F)\otimes \cA\right)_0 \ee for any complex number $\zeta\in{\mathbb{C}}$.  Denoting the vector space having as basis the set of all link diagrams on $F$ by $\mathbb{C}\cD$, we give the formula for a  morphism
\be \psi: \mathbb{C}\cD\rightarrow \mathbb{C}\cD\otimes \cA \ee that descends  to $\phi$. The extension and the formula come from  a proof analysis of March\'{e}'s work.

The main result of this paper is the theorem relating the algebras $K_{\pm {\bf i}}^0(M)$  to the algebras  where the value of the parameter is $\pm 1$.

{\bf Theorem \ref{atnegeye}} {\em  If $M$ is a compact oriented three-manifold that has no $2$-torsion in $H_1(M;\bZ)$ then $K_{\pm {\bf i}}^0(M)$ is isomorphic to $K_{\pm 1}^0(M)$. } 

The isomorphism is canonical in the sense that it does not require additional data to determine it.  This implies that $K_{\pm{\bf i}}^0(M)$  is the unreduced coordinate ring of the variety of $P\G$-characters of $\pi_1(M)$ that lift to $\G$-representations.  Our proof is based on applying the formulas  from Theorem  \ref{bige} to study a presentation of $K_{-{\bf i}}^0(M)$ coming from  a generalized Heegaard splitting of $M$.

There is related work by Sikora, \cite{S}.  If $M$ is an oriented three-manifold that can be embedded in a rational homology sphere that has no $2$-torsion in its first homology, Sikora defines an algebra structure on $K_{{\bf i}}(M)$ and then proves that it is isomorphic to $K_1(M)$.  The restriction of his algebra structure to $K_{{\bf i}}^0(M)$ agrees with the algebra structure defined here. He remarks that his theorems cannot be easily extended to all 3-manifolds since there are $3$-manifolds that are not sub-manifolds of a rational homology sphere.  We note that although indeed the algebra structure cannot be naturally defined for  $K_{\pm{\bf i}}(M)$,  the algebra structure defined in this paper  on $K_{\pm {\bf i}}^0(M)$ works for any arbitrary oriented  $3$-manifold and it is naturally defined.

The paper is structured as follows: Section \ref{charvar} contains a review of the relevant material about representation varieties and character varieties. Section \ref{skal}  recalls the definitions and needed facts about  the structure of the Kauffman bracket skein algebras of surfaces and skein modules of three-manifolds. In Section \ref{ski}  the algebra structure on $K_{\pm{\bf i}}^0(M)$ for any oriented three-manifold is introduced. In Section \ref{JM} we present the results of March\'{e} and give the formula for his isomorphism on diagrams.  Section 6 gives the proof of  the isomorphism between $K_{-{\bf i}}^0(M)$ and $K_{-1}^0(M)$ for any closed oriented three-manifold having no $2$-torsion in its first homology. We finish with a counterexample when $M$ has $2$-torsion in its first homology.
\section{ $\G$-Character Varieties} \label{charvar}
This section recalls some classical concepts of $\G$ and $P\G$ representation theory.
Let $\pi$ be a finitely generated group.
There is a naturally defined commutative  algebra \cite{LM} \be \cR^{SL_2(\mathbb{C})}(\pi)\ee  called the {\bf $\G$-representation ring of $\pi$}.  
The representation ring is constructed from the  coordinate ring of the Cartesian product of copies of $\G$, one for each generator of $\pi$, by taking the quotient  by the ideal generated by the coefficients of formal matrices obtained from instantiating the relations of $\pi$.
Representations $\rho:\pi\rightarrow SL_2(\mathbb{C})$ are in one-to-one correspondence with algebra homomorphisms \be \phi:\cR^{SL_2(\mathbb{C})}(\pi)\rightarrow \mathbb{C}.\ee

There is a right action of $PSL_2(\mathbb{C})$ on $\cR^{SL_2(\mathbb{C})}(\pi)$ coming from conjugation.  The fixed subalgebra
\be \mathcal{X}^{\G}(\pi) \ee is the {\bf universal $\G$-character ring of $\pi$}.
We say that two representations $\rho_1,\rho_2:\pi \rightarrow SL_2(\mathbb{C})$ are {\bf trace equivalent} if for every $\alpha \in \pi$, 
\be \Tr(\rho_1(\alpha))=\Tr(\rho_2(\alpha)),\ee  where $\Tr:\G\rightarrow \mathbb{C}$ is the standard trace of a matrix.
The homomorphisms from  the universal character ring to the complex numbers  are in one-to-one correspondence with trace equivalence classes of representations of $\pi$ into $SL_2(\mathbb{C})$. 

Since $\pi$ is finitely generated the algebras $\cR^{SL_2(\mathbb{C})}(\pi)$ and $ \mathcal{X}^{\G}(\pi)$ are affine over the complex numbers.  Consequently algebra morphisms from these algebras to the complex numbers are in one-to-one correspondence with the maximal ideals of the algebras via identifying a morphism with its kernel. The set of all maximal ideals of an algebra is called its {\bf maximal spectrum}.  The maximal spectrum of an affine algebra can be realized as an algebraic subset of $\bC^n$ for some $n$. Define  the {\bf $SL_2(\mathbb{C})$-representation variety}, denoted $R^{\G}(\pi)$, and {\bf $\G$-character variety, $X^{\G}(\pi)$},  to be  affine algebraic sets having the maximal spectra of $\cR^{SL_2(\mathbb{C})}(\pi)$ and
$ \mathcal{X}^{\G}(\pi) $ as their points respectively.  These sets are unique up to isomorphism of affine algebraic sets.  In general, the  representation variety (and the character variety) may not in fact be varieties, only  affine algebraic sets, but it is traditional to use the word "variety" when referring to these spaces.  

Recall 
 $\PG$ is the quotient of $\G$  by its center 
$Z(\G)=\left<\pm Id \right>$.
  A representation 
\be \rho:\pi\rightarrow \PG \ee  lifts to a representation \be \tilde{\rho}:\pi \rightarrow \G\ee if $\tilde{\rho}$ followed by the quotient map from $\G$ to $\PG$ is equal to $\rho$. \begin{definition}There is an action of $H^1(\pi;\mathbb{Z}_2)$ on $X^{\G}(\pi)$  by {\bf twisting}. If $\rho:\pi \rightarrow \G$ represents a trace equivalence class of representations  $[\rho]$ and $\alpha \in H^1(\pi;\mathbb{Z}_2)$ then
$ \alpha.[\rho]$ is the trace equivalence class of the representation that sends each $\gamma \in \pi$ to 
$ \alpha(\gamma)\rho(\gamma)$.  
\end{definition}

The fixed subalgebra of $\mathcal{X}^{\G}(\pi)$ under the action of $H^1(\pi;\mathbb{Z}_2)$,
\be \label{charring} \mathcal{X}^{\PG}(\pi)=\left(\mathcal{X}^{\G}(\pi)\right)^{H^1(\pi;\mathbb{Z}_2)}\ee  is called the ring of $P\G$-characters that lift to $\G$-characters. 
The maximal spectrum of $\mathcal{X}^{\PG}(\pi)$ is denoted $X^{P\G}(\pi)$. It is called the variety of $P\G$-characters of $P\G$-representations of $\pi$ that lift to $\G$-representations.

The map 
\be\label{quotma} X^{\G}(\pi)\rightarrow X^{\PG}(\pi) \ee that takes the character of
$\tilde{\rho}:\pi\rightarrow \G$ to the character of $\rho:\pi\rightarrow \PG$ is a regular branched cover of its image.  The group of deck transformations is   \be  H^1(\pi;\mathbb{Z}_2)=\{\alpha:\pi\rightarrow \langle\pm 1\rangle\}, \ee where the maps $\alpha$ are homomorphisms to the multiplicative group $\langle\pm 1\rangle$. 


\section{Kauffman bracket skein module}\label{skal}

In this section the  definition and properties of skein modules are recalled,  and the relationship between a skein module of a compact oriented $3$-manifold and the unreduced coordinate ring of $PSL_2(\mathbb{C})$-characters of its fundamental group that lift to $\G$-characters is explicated.

\def\KzM{K_\zeta(M)}
\def\KzuM{K_\zeta^u(M)}

Let $M$ be a  connected oriented three-manifold, and $\zeta\in\mathbb{C}-\{0\}$.  A {\bf framed link} in $M$ is a collection of disjoint annuli embedded in $M$.  The {\bf Kauffman bracket skein module of $M$} is the quotient of the complex vector space with basis the isotopy classes of framed links in $M$  by the subvector space spanned by the Kauffman bracket skein relations:  
\begin{align}
\lcr - \zeta \zer- \zeta ^{-1}\ift ,  \label{KBSR}\\
\bigcirc \sqcup L  +  (\zeta^2 +\zeta^{-2}) L\label{KBSR2}.
\end{align} 
The relations are linear combinations of framed links that are identical outside the ball where the diagrams reside. Assume that the framed link intersects the ball in strips that are parallel to the arcs in the diagram, so that if the two strips come from the same component, the same side of the annulus is up.

 Since the relations \eqref{KBSR} and \eqref{KBSR2} preserve the $\BZ_2$-homology class, there is a $H_1(M;\BZ_2)$-grading of $\KzM$:
\be \label{homgrading}
\KzM = \bigoplus_{u\in H_1(M;\BZ_2)} \KzuM,
\ee
where $\KzuM$ is the subspace spanned by  framed links $\al\in M$ such that the core of $\al$  represents $u\in H_1(M;\BZ_2)$. \begin{definition}\label{Kzero}There is an action of $H^1(M;\mathbb{Z}_2)$ on $K_\zeta(M)$ described as follows. Given $c:H_1(M;\mathbb{Z}_2)\rightarrow \langle\pm 1\rangle$, it acts as multiplication by $c(u)$ on the subspace $K^u_\zeta(M)$. The subalgebra fixed by this action is $K^0_\zeta(M)$.
\end{definition}

In general, $K_\zeta(M)$ is not an algebra.  However, in some cases  it is. 

\begin{fact} \label{three} For any oriented $3$-manifold $M$ the module $K_{\pm 1}(M)$ has an algebra structure with the product coming from perturbing framed links so that they are disjoint and then taking their  union \cite{Bu}. 

When $M=F\times [0,1]$ for a surface $F$ then $K_\zeta(M)$ has an algebra structure that comes from placing one link above the other in the direction of the interval and extending bilinearly to skeins.

 When $\zeta =\pm 1$ the stacking product and disjoint union products  on $K_\zeta(F\times [0,1])$ 
 coincide.  
\end{fact}
 In \cite{PS}, it is proved that the algebra $K_{- 1}(M)$ is naturally isomorphic to the universal $SL_2(\mathbb{C})$-character  ring  of $\pi_1(M)$.
Given a choice of a spin structure on $M$ there is an isomorphism $K_{\zeta}(M)\rightarrow K_{-\zeta}(M)$, \cite{Ba}. Consequently the algebras $K_{-1}(M)$ and $K_1(M)$ are both isomorphic to $\mathcal{X}^{\G}(\pi_1(M))$.

The action of $H^1(M;\mathbb{Z}_2)$ on $K_{\pm 1}(M)$ intertwines with the action of $H^1(M;\mathbb{Z}_2)$ on  $\mathcal{X}^{\G}(M)$.

\begin{theorem}\label{uncharcor} The ring $K_{\pm 1}^0(M)$ is the unreduced coordinate ring of the character variety of $PSL_2(\mathbb{C})$-representations of $\pi_1(M)$ that lift to representations into $\G$. \end{theorem}
 
 \proof The action of $H^1(M;\mathbb{Z}_2)$ on $K_{\pm 1}(M)$ preserves grading by $H^1(M;\mathbb{Z}_2)$.
Therefore $K_{\pm 1}^0(M)$ is isomorphic to the fixed ring under this action, which is $\mathcal{X}^{P\G}(\pi_1(M))$ as in Equation (\ref{charring}). \qed
 
 By contrast, the disjoint union does not yield an algebra structure on a module $K_{\pm{\bf i}}(M)$  since a crossing change in a framed link changes the skein of the link by a sign.

 \section{Skein Algebras  at $\pm{\bf i}$.}\label{ski}
 Although the modules $K_{\pm{\bf i}}(M)$ are not generically algebras, there is  an algebra structure on $K^0_{\pm{\bf i}}(M)$.

Given a three-manifold $M$ denote the intersection pairing on $\bZ_2$-homology by
\be \bullet: H_1(M;\mathbb{Z}_2)\otimes H_2(M;\mathbb{Z}_2)\rightarrow \mathbb{Z}_2. \ee 
Define
\be R(M)=\langle x\in H_1(M;\mathbb{Z}_2)\ |\  \forall y \in H_2(M;\mathbb{Z}_2)\  x\bullet y=0\rangle. \ee 
By Poincar'e duality $R(M)$ is the image of $H_1(\partial M;\mathbb{Z}_2)$ by the inclusion $\partial M\subset M$.
If $M$ is closed the intersection pairing
is nondegenerate, that is,  $R(M)=\langle 0\rangle $ for closed $3$-manifolds.  
If $M$ is not closed then $R(M)$ can be nontrivial. In particular, if $F$ is a  surface then $R(F\times [0,1])=H_1(F\times [0,1];\mathbb{Z}_2)$.   Any handlebody can be seen as a cylinder over a surface, consequently when $H$ is a handlebody $R(H)=H_1(H;\mathbb{Z}_2)$. 

If $L$  and $L'$ are disjoint links such that $L$ represents $0$ in $H_1(M;\BZ_2)$ and $L'$ represents an element of $R(M)$,  define the $\BZ_2$-linking number $\lk_2(L,L') \in \langle 0,1\rangle $ as follows: If  $U$ is any surface in $M$ transverse to $L'$ such that $\partial U=L$ the $\mathbb{Z}_2$-linking number is
\be 
\lk_2(L,L')= | U \cap L'| \pmod 2. 
\ee
 If $L$ and $L'$ both represent $0$ in $H_1(M;\mathbb{Z}_2)$ then $\lk_2(L,L')=\lk_2(L',L)$. Define 
  \be \label{kr} K^r_{\zeta}(M) = \bigoplus_{x\in R(M)} K^x_{\zeta}(M).\ee

 \begin{prop} \label{algstru}
 Suppose $\al\in K^0_{\pm{\bf i}}(M)$ and $\al'\in K^r_{\pm{\bf i}}(M)$   are  represented  by disjoint framed links $L$ and $L'$ respectively.  The skein  of $(-1)^{\lk_2(L,L')} L \sqcup L'$ in $K_{\pm{\bf i}}^r(M)$  is independent of the choice of $L$ and $L'$.  Consequently the assignment
  \be(\alpha,\beta) \mapsto [(-1)^{\lk_2(L,L')} L \sqcup L']  \ee
  extends bilinearly to a well-defined product
   \be\label{manal} \mu: K_{\pm{\bf i}}^0(M)\otimes K_{\pm{\bf i}}^r(M)\rightarrow K_{\pm{\bf i}}^r(M). \ee 
The restriction of $\mu$ to $K_{\pm{\bf i}}^0(M)\otimes K_{\pm{\bf i}}^0(M)$ makes $K_{\pm{\bf i}}^0(M)$ into a commutative algebra. In addition,  $\mu$ gives  $K_{\pm{\bf i}}^r(M)$ the structure of a module over $K^0_{\pm{\bf i}}(M)$. 
  \qed
   \end{prop}     
   If there is a homeomorphism between $M$ and $F\times [0,1]$, the product structure allows us to define a second product on $K^0_{\pm{\bf i}}(M)$ from stacking (discussed  in Fact \ref{three}). The two products coincide. Further, if $N\subset M$ then the map induced by the inclusion of $N$ into $M$ induces a morphisms of algebras
   \be \inc:K_{\pm{\bf{i}}}^0(N)\rightarrow K_{\pm{\bf i}}^0(M).\ee

 \section{Understanding a theorem  of March\'{e}}\label{JM} 
  This section  presents  the isomorphism constructed in \cite{Ma}, relating skein algebras of surfaces at $\zeta$ equal to $-{\bf i}$ and $-1$. We  prove an extension of the theorem that gives an isomorphism between skein algebras of surfaces when the value of the parameter $\zeta$ is twisted by ${\bf i}$
 and derive a formula for computing the value of the isomorphism on the skein induced by a diagram.

Recall that a finite type surface is a closed oriented surface with a finite number of points removed from it. In the whole section we are working with three-manifolds that are a product of a finite type surface $F$ with an interval. Since the algebra structure depends on how a $3$-manifold is presented as such a product, we emphasize it by shortening the notation to $K_\zeta(F\times [0,1])= K_\zeta(F)$.

In order to define March\'{e}'s isomorphism we need to introduce an algebra $\cA$. 
To start, given  a finite type surface $F$  consider the  vector space $\overline{\cA}$ over 
the field of complex numbers that has  basis $\{[\gamma]\in H_1(F;\mathbb{Z})\}$, where the closed braces denote the homology class represented by   the cycle $\gamma$.

The  pairing  
\be \cdot :H_1(F;\bZ)\otimes H_1(F;\bZ)\rightarrow \mathbb{Z},\ ([\gamma],[\eta])\mapsto \gamma \cdot \eta,\ee 
given by the  algebraic intersection number,  is bilinear and antisymmetric.
 There exists an associative product on $\overline{\cA}$ given by the distributive extension of the formula
 \be\label{productinA} [\gamma][\eta]={\bf i}^{-\gamma\cdot \eta}[\gamma+\eta].\ee
\begin{definition}
Define the algebra $\cA$ to be the quotient of $\overline{\cA}$ by the two sided ideal generated by the relations
  \be [\gamma]^2=[2\gamma]=1\ee for all $[\gamma]\in H_1(F;\mathbb{Z})$. 
\end{definition}
The algebra $\cA$ is graded by $H_1(F;\mathbb{Z}_2)$.   Any choice of representatives of the elements of $H_1(F;\mathbb{Z}_2)$ forms a basis of $\cA$.

 If $[\gamma],[\eta]\in H_1(F;\mathbb{Z})$ represent the same element of $H_1(F;\mathbb{Z}_2)$ then
\be [\gamma]-[\eta]=[2\beta]\ee for some cycle $\beta\in H_1(F;\mathbb{Z})$.   From the formula (\ref{productinA}) 
\be [\gamma]=[\eta+2\beta]={\bf i}^{\eta\cdot 2\beta}[\eta][\beta]^2={\bf i}^{\eta\cdot 2\beta}[\eta].\ee

In the algebra  $\cA$  we then have $[\eta]=\pm [\gamma]$ since $\eta\cdot 2\beta$ is  even. Using the facts that $\eta \cdot \eta=0$ and $\gamma=\eta+ 2\beta$ the formula becomes
\be \label{seven} [\gamma]={\bf i}^{\eta\cdot \gamma}[\eta].\ee

The tensor product
\be K_{-1}(F)\otimes_{\mathbb{C}}\cA \ee is an algebra under the coordinate-wise product. 
Since both  algebras $K_{-1}(F)$ and $\cA$ are graded by $H_1(F;\mathbb{Z}_2)$ their tensor product is graded by the sum of the gradings on the two factors.
March\'{e} defines the  {\bf diagonal subalgebra}   to be the $0$-graded part of $ K_{-1}(F)\otimes_{\mathbb{C}}\cA$.
Note that \be \left( K_{-1}(F)\otimes_{\mathbb{C}}\cA\right)_0 \ee 
 is spanned by tensors where both factors represent the same element of $H_1(F;\mathbb{Z}_2)$.

 Let $\bC\cD$ denote the vector space with basis equal to the set of isotopy classes of link diagrams on $F$. A diagram is {\bf simple} if it has no crossings and  no component of the diagram bounds a disk in  $F$. Let $\cS$ denote the set of isotopy classes of simple diagrams on $F$.  Let $\bC\cS$ be the vector space having $\cS$ as basis.  For any value of $\zeta\neq 0$,   the underlying vector space  of $K_\zeta(F)$ is 
$\bC\cS$. 
 There is a map
$ \cS \rightarrow \cA $  defined as follows.
 Given  $\alpha \in \cS$ orient the components of $\alpha$ to get a $1$-cycle $\overline{\alpha}$, then
\be \alpha \rightarrow [\overline{\alpha}].\ee   
This map is well defined since any two ways of orienting the components of $\alpha$ differ by a cycle that has zero intersection with both.

\begin{theorem}[\cite{Ma}]\label{Marorig} Let $F$ be a finite type surface  and  $\mathbb{C}\cS$ the vector space with basis given by the set of simple diagrams on $F$ up to isotopy. For any $\alpha\in \cS$ let  $n(\alpha)$  denote the number of components of $\alpha$ and $\ov{\alpha}$ denote  an arbitrary choice of orientations of the components of $\alpha$. 
The map $\phi:\bC\cS\rightarrow \bC\cS\otimes \cA$ defined by
\be \label{marc} \phi(\alpha)=(-1)^{n(\alpha)}\alpha\otimes [\ov{\alpha}]\ee
 yields a well defined isomorphisms of algebras
  \be \phi:K_{-{\bf i}}(F)\rightarrow  \left( K_{-1}(F)\otimes_{\mathbb{C}}\cA\right)_0.\ee \end{theorem}  \qed

\begin{figure}[H]\begin{picture}(91,57) \includegraphics{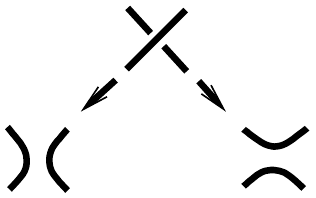}\put(-91,-5){$+1$}\put(-24,-5){$-1$} \end{picture} \caption{Positive and Negative Smoothings}\label{signsstates} \end{figure}

  We want to extend the formula from Theorem \ref {Marorig} to a formula that has linear combinations of arbitrary  link diagrams as its domain. 
If $D$ is a link diagram, a {\bf state} of $D$ is a choice of smoothing for each crossing.  Each smoothing is positive or negative as shown in Figure \ref{signsstates}.

The sum of $\pm 1$ over all the smoothings of a state $s$, where the positive smoothings contribute  $+1$ and the negative smoothings $-1$, is denoted $c(s)$. The number of components of $s$ that bound disks is $t(s)$. Let $s'$ denote the simple diagram obtained from $s$ by deleting all the trivial components. The Kauffman bracket is the map
\be\label{Kbr} \langle \ \rangle_\zeta : \bC\cD\rightarrow \bC\cS \ee given by
\be \langle D\rangle_\zeta=\sum_{s \ \mathrm{state}}\zeta^{c(s)}(-\zeta^2-\zeta^{-2})^{t(s)} s' .\ee
The Kauffman bracket cannot see isotopies  and the Reideimeister II and III moves. It gives rise to a convenient alternative definition of the the Kauffman bracket skein algebra of an oriented surface.

\begin{rem}Given an oriented surface $F$   define the algebra $K_\zeta(F)$ to be the vector space $\bC\cS$  with product defined by  the bilinear extension of the operation given by
  stacking  two simple diagrams and then applying the bracket $\langle\ \rangle_\zeta$ to the result.\end{rem}

\begin{figure}[H]\begin{picture}(91,57)\includegraphics{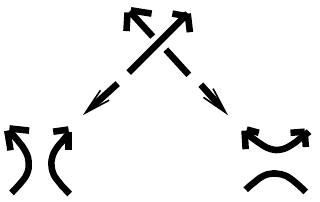}\end{picture}\caption{The Seifert smoothing is on the left}\label{seifertsmoo}\end{figure}

Given a diagram $D$  let $\ov{D}$ denote an oriented diagram coming from choosing an orientation for each component of $D$. Further denote by $n(D)$  the number of components of $D$,  $cr(D)$ the number of crossings  and $w(\ov{D})$ the writhe of $\ov{D}$. 
Recall that the Seifert smoothing of  an oriented crossing is  such that the local orientations fit together after smoothing as shown in  Figure \ref{seifertsmoo}.

Given any state $s$ of $D$  let $ns(s)$ be the number of non-Seifert smoothings and $ss(s)$ be the number of Seifert smoothings that were performed to obtain $s$ from $\ov{D}$.  It is worth noting that
\be ss(s)+ns(s)=cr(D).\ee

  If the sign of a crossing is positive then the Seifert smoothing is the positive smoothing. 
Analogously, if the sign of a crossing is negative so is the Seifert smoothing, thus
\be i^{c(s)+w(\ov{D})}=(-1)^{ss(s)}.\ee

   The following is an extension of March\'{e}'s theorem that follows from his proof.
  
  \begin{theorem} \label{bige} Let $F$ be an oriented finite type surface. For any $\zeta\in \bC-\{0\}$ the mapping $\phi$ given in equation (\ref{marc}) induces an isomorphism of algebras
  \be \phi: K_{\bi \zeta}(F)\rightarrow \left( K_\zeta(F)\otimes \cA\right)_0.\ee 
  Let $\bC\cD$ be the vector space with basis equal to the set of all isotopy classes of  link diagrams on $F$ and  $\mathbb{C}\cS$ the vector space with basis given by the set of simple diagrams on $F$ up to isotopy.  Given a link diagram $D$  let $\Xi$ denote the diagram obtained from  $\ov{D}$ by applying Seifert smoothing  at each crossing. 
If $\psi:\bC\cD\rightarrow \bC\cD\otimes \cA$ is the linear extension of 
  \be \label{psi} \psi(D)=(-1)^{n(D)}i^{-w(\ov{D})}D\otimes [\Xi] \ee then the diagram
  
 \be \label{comm} \begin{CD}  \mathbb{C}\cD@>\psi >> (\mathbb{C}\cD\otimes \cA)_0 \\ @V<\ >_{{\bf i}\zeta}VV @V<\ >_{\zeta}\otimes IdVV \\ \bC\cS @>\phi>> (\bC\cS\otimes \cA)_0 \end{CD} \ee commutes.\end{theorem}

\proof  We focus on the commutativity of Diagram (\ref{comm}). Since the product in the Kauffman bracket skein algebra comes from stacking and projecting, the commutativity of the Diagram (\ref{comm}) implies that $\phi$ is an algebra morphism. The fact that it is an isomorphism follows from the fact that it takes a basis to a basis. Thus the theorem will follow once we establish that the diagram commutes.

 Suppose that $D$ is a diagram. Any state of $D$ will be denoted $s$, and the simple diagram corresponding to a state will be denoted $s'$.  
 We follow the notation  established earlier in this section: overline  for choosing orientations, $n(D)$ is the number of components of $D$, $t(s)$ is the  number of components bounding disks of the state $s$, $cr(D)$ is the number of crossings of the diagram $D$,   $w(\ov{D})$ is writhe of the oriented diagram $\ov{D}$, $c(s)$ is the sum of signs of smoothings of the state $s$, $ns(s)$ is the  number of non-Seifert smoothings of the state $s$ with respect to $\ov{D}$ ,  $ss(s)$ is the number of Seifert smoothings of the state $s$ with respect to $\ov{D}$.  
 
 To prove the theorem we need to show
 \be \left(\langle \ \rangle_\zeta\otimes Id\right)\circ \psi(D)=\phi(\langle D\rangle_{\bi \zeta}).\ee
 First we expand the equation by using the formula  (\ref{marc}) for $\phi$, formula (\ref{psi}) for $\psi$, and (\ref{Kbr}) for the bracket with variables $\zeta$ and ${\bf i}\zeta$ respectively. 
 \be (-1)^{n(D)}\bi^{-w(\ov{D})}\sum_s\zeta^{c(s)}(-\zeta^2-\zeta^{-2})^{t(s)}s'\otimes[\Xi]=\sum_s(-1)^{n(s')}(\bi \zeta)^{c(s)}(\zeta^2+\zeta^{-2})^{t(s)}s'\otimes [\ov{s}].\ee
 This equality holds if it is true term by term, so it is enough to check that for any state $s$,
 \be (-1)^{n(D)}\bi^{-w(\ov{D})}\zeta^{c(s)}(-\zeta^2-\zeta^{-2})^{t(s)}s'\otimes[\Xi]=(-1)^{n(s')}(\bi \zeta)^{c(s)}(\zeta^2+\zeta^{-2})^{t(s)}s'\otimes [\ov{s}].\ee
 Using the fact that $[\ov{s}]=(-1)^{\frac{1}{2}\Xi\cdot \ov{s}}[\Xi] $ we see that this is equivalent to checking that
 \be (-1)^{n(D)}\bi^{-w(\ov{D})}\zeta^{c(s)}(-\zeta^2-\zeta^{-2})^{t(s)}s'\otimes[\Xi]=(-1)^{n(s')}(\bi \zeta)^{c(s)}(\zeta^2+\zeta^{-2})^{t(s)}(-1)^{\frac{1}{2}\Xi\cdot\ov{s}}s'\otimes [\Xi].\ee
 Since the images of simple diagrams under $\phi$ form a basis of $(\bC\cS\otimes \cA)_0$ it is enough to see that the coefficients on both sides agree. That means we are checking if 
 \be (-1)^{n(D)}\bi^{-w(\ov{D})}\zeta^{c(s)}(-\zeta^2-\zeta^{-2})^{t(s)}=(-1)^{n(s')}(\bi \zeta)^{c(s)}(\zeta^2+\zeta^{-2})^{t(s)}(-1)^{\frac{1}{2}\Xi\cdot\ov{s}}\ee is true.
 If $(\zeta^2+\zeta^{-2})^2 \neq 0$  we divide both sides by $(\zeta^2+\zeta^{-2})^t\zeta^{c(s)}$. If not then we just divide by $\zeta^{c(s)}$. In any case we are left to check whether 
 \be (-1)^{n(D)}\bi^{-w(\ov{D})}(-1)^{t(s)}=(-1)^{n(s')}\bi ^{c(s)}(-1)^{\frac{1}{2}\Xi\cdot\ov{s}}.\ee
 Moving $\bi^{-w(\ov{D})}$ to the right hand side yields
 \be (-1)^{n(D)}(-1)^{t(s)}=(-1)^{n(s')}\bi ^{c(s)+w(\ov{D})}(-1)^{\frac{1}{2}\Xi\cdot\ov{s}}.\ee
Using the fact that $(\bi )^{c(s)+w(\ov{D})}=(-1)^{ss(s)}$ and moving everything to one side we get
\be   (-1)^{n(D)+t(s)+n(s')+ss(s)+\frac{1}{2}\Xi\cdot \ov{s}}=1.\ee

Note $t(s)+n(s')=n(s)$, so we are really asking if 
\be  \label{almost} (-1)^{n(D)+n(s)+ss(s)+\frac{1}{2}\Xi\cdot \ov{s}}=1.\ee

In order to show that this quantity is equal to $1$ we use the following Lemma.

\begin{lemma} (Lemma 3.3. in \cite{Ma}) Let $G$ be a finite graph such that in the neighborhood of any vertex, the
edges incoming to that
 vertex have a cyclic order. We decompose the edges of $G$ in two
parts: \be E(G) = E_h \cup E_m.\ee  The edges of $E_h$ will be said to be of type handle whereas
the edges of $E_m$ will be of type M\"{o}bius.
We construct a surface $S$ from these data in the following way: take a family of
oriented discs parametrized by vertices of $G$. For all edges, we attach a band to the
corresponding discs such that the cyclic orientation of the vertices is respected. The
band should respect the orientations of the discs if the edge has type handle and should
not respect them if the edge has type M\"{o}bius.
Orient the boundary of $S$ in an arbitrary way. Let $n$ be the number of boundary
components and $m$ be the number of M\"{o}bius bands whose sides are oriented in the
same direction. Then the following formula holds:
\be n + m + \chi(S) = 0
\pmod{2},\ee where $\chi(S)$ denotes the Euler characteristic of $S$.
\end{lemma}\qed

Given the diagram $D$ form a surface $S$    as follows.  Choose an oriented  disk for each component of $\ov{D}$ so that the boundary orientation coincides with the orientation of the component.   Mark the crossings of the diagram on the  boundaries of the disks.  Add a band (rectangle)  for every crossing, where the band preserves orientation if $s$ was obtained from $\ov{D}$ by a Seifert smoothing at this crossing, and reverses it otherwise.  The boundary components of the resulting surface correspond to the components of $s$. Orient them corresponding to the orientation on $\ov{s}$.

\begin{figure}[H] \includegraphics{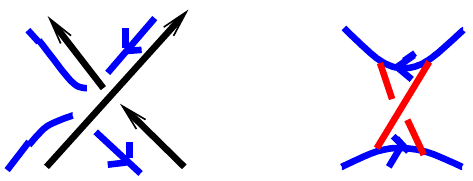}  \caption{A non-Seifert smoothing that contributes zero to $\Xi\cdot \ov{s}$.}\label{parallel} \end{figure}

We need to understand non-Seifert smoothings that  contribute zero  to $\frac{1}{2}\Xi\cdot\ov{s}$. Any state can be obtained from the diagram $\ov{D}$ by $1$-surgery. One can think of this as gluing a rectangle to the diagram $\ov{D}$ near the crossing along two of its opposite  sides, and replacing the two arcs of $\ov{D}$ that  the rectangle was glued along by the other two edges of the rectangle.  In Figure \ref{parallel} we show on the left a non-Seifert smoothing that contributes $0$ to $\Xi\cdot \ov{s}$. On the right we show a rectangle that we  have glued to the diagram $\ov{D}$ to get the state $s$.  The two blue arcs at the top and the bottom are the free edges of the rectangle and the red arcs are where it is attached to $\ov{D}$. The blue arcs have been oriented to agree with the orientaiton on $\ov{s}$. Notice that when you flatten the rectangle the arrows on the blue sides  point in the same direction. In fact,  a non-Seifert smoothing contributes zero to the algebraic intersection number $\frac{1}{2}\Xi\cdot \ov{s}$ only  when the two sides of the rectangle have the same direction in the orientation inherited from $\ov{s}$. 

From the discussion above,  $\frac{1}{2}\Xi\cdot \ov{s}$ is the number of non-Seifert smoothings of $\ov{s}$ where the two free  arcs  of the surgery rectangle are not oriented in the same direction.  Recalling  $cr(D)$ is the number of crossings of $D$ and following the notation for $m$ from the lemma, note  that 
\be cr(D)+ss(\ov{s})+\frac{1}{2}\Xi\cdot \ov{s}=m\ee is the number of non-Seifert smoothings where in the orientation of $\ov{s}$ the free arcs of the surgery rectangle have the same direction. We add   and subtract $cr(D)$ to the exponent on left hand side of Equation (\ref{almost}) to get
\be\label{bigenchilada} (-1)^{(n(D)-cr(D))+n(s)+m}.\ee

The Euler characteristic of the surface  $S$ is $n(D)-cr(D)$, and $m$ is the number of M\"{o}bius bands whose sides have the same direction. By Lemma 3.3 of \cite{Ma},  the  number in Equation (\ref{bigenchilada}) is equal to $1$. \qed

\begin{cor} \label{Mariso} The map $\psi:\bC\cD\rightarrow \bC\cD\otimes \cA$ from Equation (\ref{psi}) descends to give the isomorphism $\phi$ of March\'{e}  between the algebras  $K_{-{\bf i}}(F)$ and $\left( K_{-1}(F)\otimes_{\mathbb{C}}\cA\right)_0.$
Its restriction yields an isomorphism between 
 $K_{-{\bf i}}^0(F)$ and $K_{-1}^0(F)$. Therefore \be K_{\pm{\bf i}}^0(F)=\mathcal{X}^{\PG}(\pi_1(F)).\ee \end{cor} \qed
 
 \begin{rem}
It should be noted that although the map $\psi$ descends to skein algebras, it does not preserve skein relations for crossing involving the same component of a diagram.
 \end{rem}



 \section{ The isomorphism between $K_{-{\bf i}}^0(M)$ and $K_{-1}^0(M)$.}
 In this section we use the isomorphism $\phi$ from Corollary \ref{Mariso} to prove the analogous result for  $3$-manifolds that are not necessarily a cylinder over a surface. Given a compact oriented $3$-manifold $M$ recall that $K_{-{\bf i}}^0(M)$ has an algebra structure defined by  the restriction of \eqref{manal}.

 \begin{theorem}\label{atnegeye}
 Let $M$ be a compact oriented $3$-manifold such that $H_1(M;\mathbb{Z})$ has no $2$-torsion. The skein algebras $K_{\pm{\bf i}}^0(M)$ and $K_{\pm1}^0(M)$ are isomorphic.
 
\end{theorem}

We will focus on $K_{-{\bf i}}^0(M)$ and $K_{-1}(M)$. The proof in the other case follows by substitution.

 We prove that given a generalized Heegaard surface $F$ for the $3$-manifold $M$,  the restriction of $\phi: K_{-\bi }(F)\rightarrow \left( K_{-1}(F)\otimes \cA\right)_0$ to 
 \be \phi:K_{- {\bf i}}^0(F)\rightarrow \left(K_{-1}^0(F)\otimes \cA\right)_0\ee  descends  to an isomorphism 
 \be \phi:K^0_{-\bf i}(M)\rightarrow (K^0_{-1}(M)\otimes \cA)_0=K^0_{-1}(M).\ee There are examples of three-manifolds with $2$-torsion in $H_1(M;\mathbb{Z})$ for which this isomorphism does not descend. An example of  how this happens is shown in Subsection \ref{counterexample}.

The proof of {\bf  Theorem \ref{atnegeye}}  uses the concept of  {\bf handle-slides}. We begin by explaining what these are.

Let $M$ be a compact oriented three-manifold and  let $F$ be be a {\bf generalized  Heegaard surface} for $M$. This means that $M$ is homeomorphic to the result of adding $2$-handles to both sides of $F\times [0,1]$ and then perhaps capping off some sphere boundary components with balls.  Let $f:M\rightarrow [0,3]$ be a self indexing Morse function on $M$, that is  all its critical points lie in the interior of $M$, and $f^{-1}(i)$ is the set of critical points of index $i$. The existence of such a function is guaranteed by standard results about Morse functions.  The surface $f^{-1}(3/2)$ is a generalized Heegaard surface. A collar of the generalized Heegaard surface is given by $f^{-1}([5/4,7/4])=f^{-1}(3/2)\times [0,1] $ .  The $2$-handles are fattened up versions of part of ascending manifolds of $1$-handles that lie below $f^{-1}([5/4,7/4])$ and the fattened up versions of the parts of descending manifolds of index $2$-critical points that lie above $f^{-1}([5/4,7/4])$.  The balls come from neighborhoods of the index $0$ and index $3$ critical points.

Hence $M=H_1\cup H_2$ where $H_1\cap H_2=F$ and $H_1$ and $H_2$ are the results of adding $2$-handles to different sides of a collar of $F$ and maybe capping off some sphere boundary components.  This can be visualized as a diagram on the surface $F$, which  consists of  a collection of disjoint red curves $R$ that are  attaching curves for the $2$-handles on one side and a  collection of disjoint blue curves $B$ that are  attaching curves on the other side. Imagine the red curves as lying over and the blue curves as lying under the surface $F$ in a collar of  $F$. We can ignore the balls used to cap sphere boundary components.
 
 Let $\mathbb{C}\mathcal{L}(M)$ denote the vector space whose basis is the set of isotopy classes of framed links in $M$. There is a surjective linear map
   \be \label{inc}  \inc:\mathbb{C}\mathcal{L}(F)\rightarrow \mathbb{C}\mathcal{L}(M)\ee
   induced by the inclusion $F\times [0,1]\rightarrow M$ coming from identifying a collar of $F$ in $M$ with $F\times [0,1]$.

   Given a framed link $L$ representing an element of ${\cL(F)}$, a {\bf handle-slide} is the difference between $L$ and a band sum of $L$ with one of the curves in $R$ or $B$. 
    The following is a consequence of the theory of singularities of smooth mappings
  \begin{fact} 
   The  kernel of $\operatorname{inc}:\mathbb{C}\mathcal{L}(F)\rightarrow \mathbb{C}\mathcal{L}(M)$ is spanned by  handle-slides. 
 \end{fact}

   \begin{figure}[H]\begin{center} \scalebox{.5}{\includegraphics{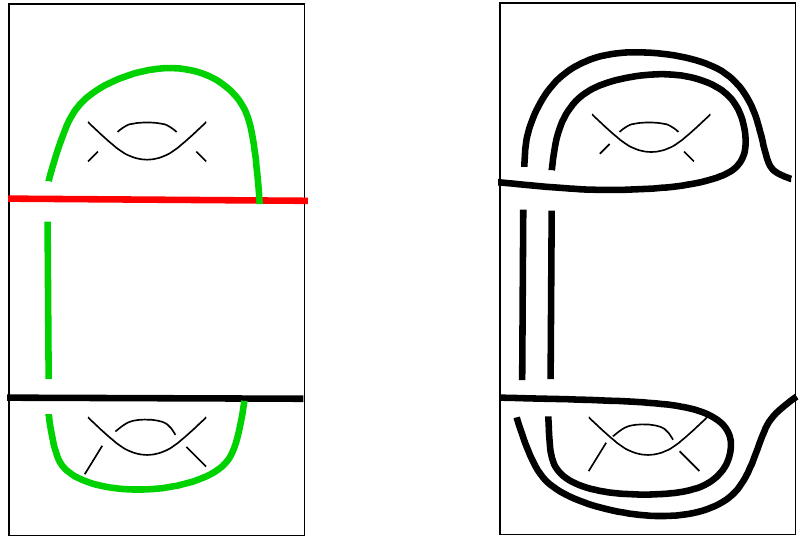}} \end{center} \caption{Elementary handle-slide}\label{elhasli} \end{figure}

    We will view  blackboard framed links in a collar of $F$ in $M$ as link diagrams on the surface $F$. In order to prove that  the  skein algebras $K_{\pm{\bf i}}^0(M)$ and $K_{-1}^0(M)$ are isomorphic we will use the map 
  $\psi$ from  Equation (\ref{psi}) defined on the set of isotopy classes of link diagrams on $F$, denoted $\mathbb{C}\cD$ .
 Hence we will interpret the map of Equation (\ref{inc}) as
   \be \inc:\mathbb{C}\cD \rightarrow \cL(M)\ee
   and we want to describe handle slides in terms of diagrams.
 
 In the language of diagrams on $F$, an {\bf elementary handle-slide} is built as follows. Let $s$ be a diagram and let $a$ be an arc in $F$ with one endpoint in $s$ and another endpoint on simple closed curve $c$ that belongs to  $R$ or  $B$. The diagram of the arc $a$ may  have over- or under- crossings with $s$ and $c$   in its interior.  In the left side of Figure \ref{elhasli} the  green arc $a$ joins one of the red curves with the diagram $s$ depicted in black.
 Replace $a$ by an embedding of a strip $N(a)$ that intersects each of $s$ and $c$  in a single arc and passes over  or  under $s$ at all the points of intersection of the arc with $s$ in its interior.   The strip always passes over  (or under) $c$ at each point of intersection of $a$ with $c$, depending on whether $c$ is in $B$ (or $R$ respectively). 
The   strip $N(a)$ that  replaces $a$  can be chosen to be blackboard framed this can be seen by making the strip arbitrarily short to start with, but then stretching it out so that it becomes a blackboard framed arc on $F$.   Also all crossings that involve only the diagram can be resolved, so we can assume that we are only working with simple diagrams \cite{B2}.
  Form a new diagram, $D$, by taking the union of $s$ with the boundary of the strip $N(a)$ and removing the arcs where the $N(a)\cap s$. The diagram $D$ is pictured on the right in Figure \ref{elhasli}. 
   The difference
   \be s-D \ee in the vector space $\mathbb{C}\mathcal{D}$ is an elementary handle-slide.

   A {\bf compound handle-slide} consists of performing several elementary handle-slides in succession, with their respective arcs $a_i$  diagrammatically missing one another by passing over or under. For subsequent slides, the arc $a_j$ does not necessarily have to pass over (resp. under) the blue (resp. red)  curves, in  case these were involved in one of the previous slides. A compound slide that involves multiple slides over the same $2$-handle is realized  by having parallel copies of each of the attaching curves. A compound handle-slide is the difference $L-L'\in \mathbb{C}\mathcal{D}(F)$ where $L$ is the starting diagram of a framed link and $L'$ is the result of the sequence of elementary handle-slides.

   The homology class in $H_1(F\times[0,1];\mathbb{Z}_2)$ of the compound handle slide $L-L'$ is  the sum of all the curves in $R$ and $B$, counted with multiplicities,  that are involved in the sequence of elementary handle slides comprising $L-L'$. We call this class {\bf the homology class of the compound handle slide}.

   The proof of Theorem \ref{atnegeye} uses the following lemma.
\begin{lemma}\label{saver} Let $A$ be a finitely generated free abelian group and let
  \be \omega:A\otimes A \rightarrow \mathbb{Z} \ee be an antisymmetric bilinear pairing on $A$. Suppose that $L,L'\leq A$ are subgroups with the property that the restrictions of $\omega$ to $L\otimes L$ and $L'\otimes L'$ are zero. Further assume that $A/(L+L')$ has no $2$-torsion.   For any $\alpha \in L$ and $\alpha'\in L'$, such that $\alpha+\alpha'=2\beta$ for some $\beta \in A$ the value of
  $\omega(\alpha,\alpha')$ is divisible by $4$. \end{lemma}

\proof   The bilinear form  $\omega$ descends to a pairing on $(L+L')/(L\cap L')$. The group $(L+L')/(L\cap L')$ is the direct sum of $\im(L)$ and $\im(L')$ under the quotient map,
\be\label{dirs} (L+L')/(L\cap L')=\im(L)\oplus \im(L'). \ee   Since $A/(L+L')$ has no $2$-torsion the fact that  $\alpha+\alpha'=2\beta$ for some $\beta \in A$ implies that $\beta\in L+L'$.  That means the image $\ov{\alpha+\alpha'}$  of $\alpha+\alpha'$   is divisible by $2$ in  $(L+L')/(L\cap L')$. There are $\gamma\in L$ and $\delta\in L'$ such that their images $\ov{\gamma}$ and $\ov{\delta}$ in $(L+L')/(L\cap L')$ have the property that
\be\ov{\alpha}+\ov{\alpha'}=2(\ov{\gamma}+\ov{\delta}).\ee By (\ref{dirs}) we have that $\ov{\alpha}=2\ov{\gamma}$, and $\ov{\alpha'}=2\ov{\delta}$,    so
\be \omega(\alpha,\alpha')=\omega(\ov{\alpha},\ov{\alpha'})=\omega(2\ov{\gamma},2\ov{\delta})=4\omega(\ov{\gamma},\ov{\delta}).\ee Hence $\omega(\alpha,\alpha')$ is divisible by $4$. \qed

 We proceed with the proof of {\bf  Theorem \ref{atnegeye}}.
\proof 
As discussed above, choose a generalized Heegaard surface $F$  for the 3-manifold $M$, with a collar of $F$ in $M$ identified with $F\times [0,1]$.  Let $R\cup B$ be a complete set of attaching curves for $2$-handles for the generalized Heegaard surface with the curves in $R$ lying in $F\times \{1\}$ and the curves in $B$ lying in $F\times \{0\}$.  The Mayer-Vietoris sequence implies that 
\be H_1(M;\mathbb{Z})=H_1(F;\mathbb{Z})/(\langle R\rangle+\langle B\rangle) \ee where $\langle R\rangle$ and $\langle B\rangle $ denote the subgroups of $H_1(F;\mathbb{Z})$ spanned by oriented versions of the curves in $R$ and $B$.

Recall  the map from Theorem \ref{Mariso},
\be \phi:K^0_{-{\bf i}}(F\times I) \rightarrow \left( K^0_{-1}(F\times I)\otimes \cA\right)_0 \ee
and the map from Equation (\ref{psi})
\be \psi:\bC\mathcal {D}\rightarrow  \left(\bC\cD\otimes \cA\right)_0,\ee 
where $\mathcal {D}$ denoted the set of framed link diagrams on $F$.

  Letting $\mathbb{C}\mathcal{D}^0$ denote the isotopy classes of framed link diagrams that represent $0$ in $H_1(F\times [0,1];\mathbb{Z}_2)$ and $\mathbb{C}\mathcal{L}^0(M)$ the isotopy classes of framed links in $M$ that represent $0$ in  $H_1(M;\mathbb{Z}_2)$, the restriction of the inclusion
   \be\label{resinc} inc:\mathbb{C}\mathcal{D}^0\rightarrow \mathbb{C}\mathcal{L}^0(M)\ee is onto.  If a framed link $L$ represents $0$ in  $H_1(M;\mathbb{Z}_2)$, it  bounds a compact surface. Isotope that surface into $F\times [0,1]$. Its boundary is a link in $F\times [0,1]$ that represents $0$ in  $H_1(F\times [0,1];\mathbb{Z}_2)$ and is isotopic to $L$.   The kernel of the map in equation (\ref{resinc}) is spanned by compound handle-slides $L-L'$ such that $L-L'=0$ in $H_1(F\times[0,1];\mathbb{Z}_2)$.
   One can always choose the diagrams representing the elementary  handle-slides involved in $L-L'$ so that the strips are parallel to $F\times \{0\}$ (i.e.,  blackboard framed).
   There is a commutative diagram
   \be\label{comb} \begin{CD}  \mathbb{C}\mathcal{D}^0 @>\inc>> \mathbb{C}\mathcal{L}^0(M) @>>> 0 \\ @V<\ >_{-{\bf i}}VV @V< \ >_{-{\bf i}}VV \\ K_{-{\bf i}}^0(F) @>\inc>>  K_{-{\bf i}}^0(M) @>>>0\end{CD}. \ee Here we are using the braces to denote the Kauffman bracket. Although we defined the bracket on diagrams, it gives rise to a map on framed links by representing each framed link as a diagram and then applying the bracket.

 If $S$ is the linear span of the compound handle-slides where the starting diagram is simple, the bands are blackboard framed and the homology class corresponding to the compound handle-slide is zero in $H_1(F\times[0,1];\mathbb{Z}_2)$  then 
   \be \langle S\rangle_{-{\bf i}}= \ker\left(\inc:K_{-{\bf i}}^0(F\times[0,1])\rightarrow  K_{-{\bf i}}^0(M)\right).\ee  This  argument follows the same reasoning as  in \cite{B2}. Hence \be \label{eight} K_{-{\bf i}}^0(M)\cong K_{-{\bf i}}^0(F)/\langle S\rangle_{-{\bf i}}.\ee
   
   Similarly, there is a commutative diagram 
   \be \begin{CD}  \mathbb{C}\mathcal{D}^0 @>\inc>> \mathbb{C}\mathcal{L}^0(M) @>>> 0 \\ @V<\ >_{-1}VV @V<\ >_{-1}VV \\ K_{-1}^0(F) @>\inc>>  K_{-1}^0(M) @>>>0.\end{CD} \ee
  
   It is easy to see that  \be \langle S\rangle_{-1}= \ker\left(\inc:K_{-1}^0(F\times[0,1])\rightarrow  K_{-1}^0(M)\right).\ee
   Hence,
   \be \label{nine} K_{-1}^0(M)\cong K_{-1}^0(F)/\langle S\rangle_{-1} .\ee

   The theorem will be proved by showing that   the map from Equation (\ref{psi}), 
    \be \psi:\mathbb{C}\mathcal{D}\rightarrow  \left(\mathbb{C}\mathcal{D}\otimes \cA\right)_0,\ee 
   has the property that
   \be \psi(S)=S\otimes [0],\ee   which,  along with  Equations (\ref{eight}) and (\ref{nine}), shows that $K^0_{-{\bf i}}(M)$ and $K_{-1}^0(M)$ are isomorphic.

Let $s$ be a simple diagram on $F$ and let $D$ be the result of a compound handle-slide on $s$ whose homology class is  $0$ in $H_1(F\times[0,1];\mathbb{Z}_2)$.  Let $\ov s$ and $\ov{D}$ be  choices of orientation on $s$ and on $D$ that agree up to isotopy. Since $s$ has no crossings, $\ov s$ is the $1$-cycle corresponding to smoothing its crossings. Let $\Xi$ denote the $1$-cycle coming from smoothing the crossings of $\ov{D}$.  Since $s$ and $D$ have the same number of components, we denote
\be n=n(D)=n(s).\ee
The compound handle-slide has the  endpoints of all the arcs lying on $s$. We can choose the strips involved in the handle-slide to have blackboard framing as discussed above.  Hence in the computation of writhe of $\ov{D}$
we can ignore the contributions of the sides of the strips as they come in canceling pairs. Therefore the writhe of $\ov{D}$ is equal to the writhe of the diagram $\ov{s}\cup \ov{r} \cup \ov{b}$, where $r$ is the union of copies of curves in $R$, and $b$ is the union of copies of curves in $B$, involved in the compound handle-slide. Note that the $r$ and $b$ curves are placed above and below $s$ respectively, moreover  $r+b=0\in H_1(F;\mathbb{Z}_2)$ and the orientations on $r$ and $b$ are chosen to agree with $\ov{D}$. Letting $\omega$ denote the intersection pairing on $H_1(F;\mathbb{Z})$ we have that the writhe is
\be w(\overline{D})=\omega(\overline{r},\overline{s})+\omega(\overline{r},\overline{b})+\omega(\overline{s},\overline{b}).\ee  We are using the fact that the diagrams $s$, $b$ and $r$ are simple hence there are no self crossings. 
Since the intersection pairing is bilinear and antisymmetric, 
\be \label{mambo}  w(\overline{D})=\omega(\overline{s},\overline{b}-\overline{r})+\omega(\overline{r},\overline{b}). \ee

The fact that $\overline{b}+\overline{r}$ represents $0$ in $H_1(F;\mathbb{Z}_2)$ implies that  $\overline{b}-\overline{r}$ does also. Hence there is a $1$-cycle $f$ such that $2f$ differs from $\overline{b}-\overline{r}$ by a boundary.  Similarly, as $\overline{s}$ represents $0$ in  $H_1(F;\mathbb{Z}_2)$, there is a $1$-cycle $g$  such that $2g$ differs from $\overline{s}$ by a boundary. 
Hence 
\be \omega(\overline{s},\overline{b}-\overline{r})=\omega(2g,2f)=4\omega(g,f). \ee 
This means the first term on the right in Equation (\ref{mambo}) is divisible by $4$. 

 Let $A=H_1(F;\mathbb{Z})$, $L=\langle R\rangle $ and $L'=\langle B\rangle$. The fact that $H_1(M;\mathbb{Z})$ has no $2$-torsion allows us to apply Lemma \ref{saver} to see that $\omega(\overline{r},\overline{b})$ is divisible by $4$. Therefore \be w(\overline{D})=0\mod{4}\ee and
\be {\bf i}^{-w(\ov{D})}=1.\ee

Recall the map $\psi:\mathbb{C}\mathcal{D}\rightarrow (\mathbb{C}\mathcal{D}\otimes \mathcal{A})_0 $ from  Equation (\ref{psi}). Since $s$ has no crossings, $w(\overline{s})=0$, and
\be \psi(s-{D})=(-1)^n {(s\otimes [\ov{s}] - \bf i}^{-w(\ov{D})}D\otimes [\Xi])=(-1)^n (s\otimes [\ov{s}] - D\otimes [\Xi]).\ee Finally,
    since the cycles $\ov{s}$ and  $\Xi$ represent $0$ in $H_1(F;\mathbb{Z}_2)$, 
        \be  \psi(s-{D})=(-1)^n(s-D)\otimes [0].\ee  Hence the linear extension of $\psi$ to linear combinations of diagrams satisfies
\be \psi(S)= S \otimes [0].\ee
Therefore $\psi$  gives rise to an isomorphism
\be K_{-{\bf i}}^0(M)\rightarrow K_{-1}^0(M).\ee \qed

\subsection{Independence from the Heegaard Splitting}
In this subsection we show that the isomorphism from Theorem \ref{atnegeye} does not depend on the choice of a Heegaard splitting of a $3$-manifold.
\begin{prop} \label{independence}
The isomorphism $ \phi:K^0_{\pm \bf i}(M)\rightarrow K^0_{\pm1}(M)$  is independent of the choice of Heegaard splitting of the $3$-manifold $M$.
\end{prop}
\proof  We only prove the proposition in the case of $-{\bf i}$ and $-1$. The other case follows by substitution.

\begin{center} \begin{figure}[H] \scalebox{.5}{\includegraphics{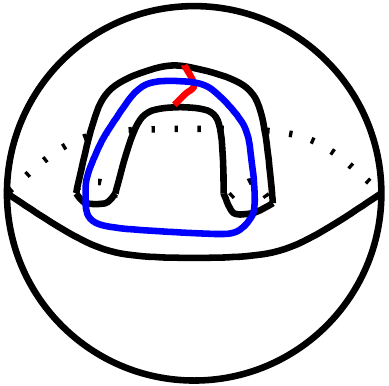}}\caption{ A trivial handle with standard meridians shown in red and blue}\label{trivh} \end{figure} \end{center}

Recall that the map $\phi$ is  a descent of the map $\phi$ defined on a Heegaard surface for $M$, 
\be\phi_F:K_{- {\bf i}}^0(F\times I)\rightarrow \left(K_{-1}^0(F\times I)\otimes \cA\right)_0.\ee
Any two Heegaard surfaces of a $3$-manifold become isotopic after a finite number of stabilizations, that is adding of a trivial handle. We show a trivial handle inside a ball in Figure \ref{trivh}. 

Suppose that the surface $F'\subset F\times [0,1]$ is  the result of adding a single trivial handle to the surface $F\times \{1/2\}\subset F\times [0,1]$.  

 A collar of $F'$ in $F$ is homeomorphic to $F'\times [0,1]$. Since $F'$ is a Heegaard surface for $F\times[0,1]$   the inclusion map $F'\times [0,1]\subset F\times [0,1]$ induces an isomorphism
 \be d:K_{-{\bf i}}(F)  \rightarrow K_{-{\bf i}}(F')/S \ee where $S$ is the submodule of $K_{-{\bf i}}(F')$ spanned by handleslides, and the map $d$ comes from isotoping framed links in $F\times [0,1]$ into the collar of $F'$ in $F$.  The restriction of this map, which we denote by the same name, is an isomorphism 
 \be d:K_{-{\bf i}}^0(F)  \rightarrow K_{-{\bf i}}^0(F')/S^0. \ee  Here $S^0$ denotes the intersection of $S$ with the $0$-graded part of $K_{-{\bf i}}(F)$.  
 Analogously, on the level of framed links, the inclusion $F'\times [0,1]\subset F\times [0,1]$ induces an isomorphism 
 \be d:K_{-1}^0(F)\rightarrow  K_{-1}^0(F')/S^0. \ee 
 If $J$ is a simple diagram on $F$ then there is an {\em almost vertical  deformation} of $J$ which is a simple diagram on $F'$ that is isotopic to $J$ and  avoids the trivial handle.     One-cycles on $F$ can be deformed to one cycles on $F'$ that miss the trivial handles. Recall that  March\'{e}'s  algebra $\cA$ was defined for a specific surface. We indicate this by using the notation $\cA(F)$.
 The deformation above gives an injective algebra morphism 
 \be \tau:\cA(F)\rightarrow \cA(F').\ee 
 
 Theorem \ref{atnegeye} implies that the map $\phi_{F'}: K_{-{\bf i}}^0(F') \rightarrow (K_{-1}^0(F')\otimes \cA(F'))_0 $ defines a map $K_{-{\bf i}}^0(F')/S^0 \rightarrow (K_{-1}^0(F')/S^0\otimes \cA(F'))_0$. We abuse notation by denoting it also by $\phi_ {F'}$.
 We need to show that the following diagram commutes.
 \be\label{twoefs} \begin{CD}  K_{-{\bf i}}^0(F) @>\phi_{F}>>  (K_{-1}^0(F)\otimes \cA(F))_0 \\ @VdVV @Vd\otimes \tau VV \\  K_{-{\bf i}}^0(F')/S^0 @>\phi_{F'}>>  (K_{-1}^0(F')/S^0\otimes \cA(F'))_0 \end{CD} \ee
 Suppose that $s$ is a simple diagram on $F$ and $s'$ is a simple diagram on $F'$ that is an almost vertical deformation of $s$. Let $\ov{s}$ and $\ov{s}'$ be oriented versions of the two diagrams such that the deformation between the two diagrams preserves orientation. Let $n$ be the number of components of each of the two diagrams. Both diagrams $\ov{s}$ and $\ov{s}'$  have writhe $0$ as they are simple.
 To see that  (\ref{twoefs})
 commutes we need to show that
 \be \phi_{F'}(d(s))=d\otimes \tau(\phi_F(s)).\ee
 Following the definitions
 \be \phi_{F'}(d(s))=\phi_{F'}(s')=(-1)^{n}s'\otimes [\ov{s}'] \ee and
 \be d\otimes \tau(\phi_F(s))=d\otimes \tau\left((-1)^ns\otimes [\ov{s}]\right)=(-1)^ns'\otimes [\ov{s}'],\ee which yields the desired result. 
 
 It is worth noting that $(K_{-1}(F)\otimes \cA(F))_0=K_{-1}^0(F)\otimes [0]$, and $(K_{-1}^0(F')/S^0\otimes \cA(F'))_0=K_{-1}^0(F')/S^0\otimes [0]=K_{-1}(F)$.

 Suppose now that $F$ is a surface yielding a generalized Heegaard splitting of an oriented three-manifold $M$ and $T\leq K_{-{\bf i}}(F)$ and (by abuse of notation) $T\leq K_{-1}(F)$ are the submodules spanned by handleslides corresponding to $M$. Let $T^0$ denote their intersection with the zero graded parts of the skein algebras.  Let $F'$ be a Heegaard surface for $M$ obtained from $F$ by adding a single trivial handle. Finally let $T'$ and ${T'}^0$ be the corresponding submodules of  $K_{-1}(F')$ (or $K_{-{\bf i}}(F')$ ).
 It is easy to see that $S^0+T^0={T'}^0$.  Hence taking quotients we get that the diagram 
 
 \be \begin{CD}  K_{-{\bf i}}^0(F)/T @>\phi_{F}>>  (K_{-1}^0(F)/T\otimes \cA(F))_0 \\ @VdVV @Vd\otimes \tau VV \\  K_{-{\bf i}}^0(F')/{T'}^0 @>\phi_{F'}>>  (K_{-1}^0(F')/{T'}^0\otimes \cA(F'))_0 \end{CD} \ee
 is commutative.  Therefore the induced diagram
 \be  \begin{CD} K_{-{\bf i}}^0(M)@>\phi_F>> K_{-1}^0(M) \\ @VdVV @VdVV \\   K_{-{\bf i}}^0(M)@>\phi_{F'}>> K_{-1}^0(M) \end{CD} \ee is commutative.

 This can be iterated through multiple stablizations. After finitely many stabilizations any two generalized Heegaard surfaces of a three-manifold become isotopic, so the isomorphism of skein modules is independent of the generalized Heegaard surface. \qed

 \begin{cor}\label{cornotor} If $M$ is a compact oriented three-manifold with no $2$-torsion in $H_1(M;\mathbb{Z})$  then the algebra $K_{\pm{\bf i}}^0(M)$ is  naturally isomorphic to the unreduced coordinate ring of characters of  $PSL_2(\mathbb{C})$-representations of $\pi_1(M)$ that lift to $SL_2(\mathbb{C})$-representations of $\pi_1(M)$.  \end{cor}

 \proof This follows from Theorem \ref{uncharcor} in the case of $K^0_{-{\bf i}}(M)$.  In the case of $K^0_{{\bf i}}(M)$
 the identification of $K_1(M)$ with $\mathcal{X}^{P\G}(\pi_1(M)) $ depends on the choice of a spin structure. However all those isomorphisms agree when restricted to $K_{1}^0(M)$. 
 \qed

 \begin{rem}In the case of ${\bf i}$ the morphism is not natural until you choose a spin structure for $M$. \end{rem}

\begin{scholium} If $M$ is a compact oriented three-manifold such that $H_1(M;\bZ)$ has no $2$-torsion, then for all $\zeta\in \bC-\{0\}$ the vector spaces $K_{{\bf i}\zeta}^0(M)$ and $K_\zeta^0(M)$ are isomorphic. \end{scholium} \qed

\subsection{ A counterexample}\label{counterexample}  As was remarked earlier, the map $\phi$ used to construct the isomorphism between  $K_{\pm {\bf i}}^0(M)$ and $K_{-1}^0(M)$ does not necessarily descend from diagrams on a generalized Heegaard surface to skeins in the $3$-manifold when there is $2$-torsion in $H_1(M;\mathbb{Z})$. The reason for this is there could be a compound handleslide where the writhe $\pmod{4}$ of the diagram after a handleslide is equivalent to $2$ while the original diagram has writhe $0$.  

\begin{figure}[H]\includegraphics{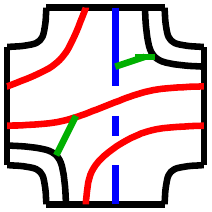} \caption{ A counterexample when $H_1(M;\mathbb{Z})$ has $2$-torsion}\label{condo} \end{figure}

In Figure \ref{condo} we show a configuration in a punctured torus that lies inside a generalized Heegaard surface. The top edge is identified with the bottom edge and the left edge is identified with the right edge. The red and blue curves are attaching curves of $2$-handles on opposite sides of the Heegaard surface. This pair of curves causes $H_1(M;\mathbb{Z})$ to have $2$-torsion. The pair of black arcs are parts of a simple diagram that pass through the punctured torus, and the two green arcs are the strips used to form the handleslide. Notice the writhe of the result of surgering the diagram along the green arcs is equivalent to $2 \pmod{4}$.

Note that the counterexample above implies that our proof of Corollary \ref{cornotor} does not work for $3$-manifolds $M$ with 
 $2$-torsion in $H_1(M;\bZ)$; it is not a counterexample to the Corollary.

\end{document}